\documentclass[a4paper]{article}
\usepackage{amsmath}
\usepackage{amssymb, amsfonts}
\usepackage{epsfig}
\usepackage{graphicx}

\begin{document}

\title{Stochastic stability of invariant measures: The 2D Euler equation}
\date{ }
\author{F. Cipriano\thanks{%
Departamento de Matem\'{a}tica, Faculdade de Ci\^{e}ncias e Tecnologia da
Universidade Nova de Lisboa and Centro de Matem\'{a}tica e Aplica\c{c}\~{o}%
es, mfsm@fct.unl.pt}, H. Ouerdiane\thanks{%
D\'{e}partment de Math\'{e}matiques, Facult\'{e} des Sciences, Universit\'{e}
de Tunis El Manar, Campus Universitaire, 1060 Tunis, Tunis,
habib.ouerdiane@fst.rnu.tn} and R. Vilela Mendes\thanks{%
CMAFCIO and IPFN, Universidade de Lisboa, C6, 1749-016 Lisboa, Portugal, \
e-mail: rvilela.mendes@gmail.com, rvmendes@fc.ul.pt,
http://label2.ist.utl.pt/vilela/}}
\maketitle

\begin{abstract}
In finite-dimensional dissipative dynamical systems, stochastic stability
provides the selection of the physical relevant measures. That this might
also apply to systems defined by partial differential equations, both
dissipative and conservative, is the inspiration for this work. As an
example the 2D Euler equation is studied. Among other results this study
suggests that the coherent structures observed in 2D hydrodynamics are
associated to configurations that maximize stochastically stable measures
uniquely determined by the boundary conditions in dynamical space.
\end{abstract}

\section{Introduction}

The main purpose of this research is to extend the notion of stochastically
stable invariant measure to dynamical systems defined by partial
differential equations, in particular to conservative systems with many
invariant measures where the notion of stochastic stability may provide a
selection criteria for the physically relevant measures. In the following
subsections and also partly in Section 2 and 3 some standard material is
formulated in a notation appropriate for the further developments. The main
original results are contained in the Section 4 and 5. The most direct
physical implication would be the interpretation of the coherent structures
observed in two-dimensional and quasi-two-dimensional fluid motion as
configurations maximizing stochastically stable invariant measures.
According to the results, the stochastically stable invariant measures would
be unique for each choice of boundary conditions in the dynamical variables.

\subsection{The physical relevance of stochastically stable invariant
measures}

For finite-dimensional systems the notions of \textit{physical measure} and 
\textit{stochastically stable measure} are closely related. Let $M$ be the
state space, $f:M\rightarrow M$ a dynamical system defined by a smooth
transformation and $\mu $ a positive Borel measure on $M$ such that 
\begin{equation}
\lim_{n\rightarrow \infty }\frac{1}{n}\sum_{j=0}^{n-1}\varphi \left(
f^{j}\left( x\right) \right) \rightarrow \int_{M}\varphi d\mu  \label{S1}
\end{equation}%
for a positive measure set $A$ of initial points $x$ and any continuous
function $\varphi :M\rightarrow \mathbb{R}$. It means that time averages of
continuous functions are given by the corresponding spatial averages
computed with respect to $\mu $, at least for a large set of initial states $%
x$. Such measure $\mu $, when it exists, is called a \textit{physical measure%
} (or Sinai-Bowen-Ruelle, SBR measure) \cite{Sinai} \cite{Bowen} \cite%
{Ruelle1} \cite{Ruelle2}.

For uniformly hyperbolic systems there is a complete theory concerning
existence and uniqueness of physical measures and partial results for
non-uniformly hyperbolic and partially hyperbolic systems \cite{Katok} \cite%
{Barreira}.

Consider now the stochastic process $f_{\varepsilon }$ obtained by adding a
small random noise to the deterministic system $f$. Under very general
conditions, there exists a stationary probability measure $\mu _{\varepsilon
}$ such that, almost surely,%
\begin{equation}
\lim_{n\rightarrow \infty }\frac{1}{n}\sum_{j=0}^{n-1}\varphi \left(
f_{\varepsilon }^{j}\left( x\right) \right) \rightarrow \int \varphi d\mu
_{\varepsilon }  \label{S2}
\end{equation}%
\textit{Stochastic stability of the }$\mu $\textit{\ measure} means that $%
\mu _{\varepsilon }$ converges to the physical measure $\mu $ when the noise
level $\varepsilon $ goes to zero. There is stochastic stability for
uniformly hyperbolic maps, for Lorenz strange attractors, H\'{e}non strange
attractors and also general results for partially hyperbolic systems \cite%
{Kifer1} \cite{Young} \cite{Kifer2} \cite{Viana1} \cite{Viana2} \cite{Alves}%
. Existence and uniqueness of the invariant measure $\mu _{\varepsilon }$
under general conditions provides a powerful tool to obtain the relevant
physical measure of the dynamical system $f$, by randomly perturbing it and
then letting the noise level $\varepsilon \rightarrow 0$.

In the past, stochastic stability of the physical measures has been
considered mostly relevant for dissipative systems or for Hamiltonian
systems with small dissipative perturbations. That the same notion might
also be useful for strictly conservative systems follows from our results,
with the choice of boundary conditions in the dynamical space leading to
uniqueness of the stochastically stable measure.

\subsection{The 2D Euler equation and persistent large-scale structures in
(quasi) two dimensional fluid motion}

For definiteness, our study concentrates on the stability of invariant
measures for the 2D Euler equation, an issue of current physical interest
for the understanding of geophysical phenomena \cite{Bouchet1} \cite%
{Bouchet2}. A striking feature of (quasi) two-dimensional turbulent fluid
motion \cite{Bouchet1} is the emergence of large scale structures which
persist for long time intervals. Another feature is the relaxation of the
flow to a small number of patterns, as if they were attractors of the
dynamics, a feature not to be expected in conservative or small dissipation
systems. This last feature, is also contrary to the idea that viscosity is
required to explain irreversibility in turbulent flows. These phenomena
should hopefully be explained by the 2D Euler equation or by its
quasi-geostrophic variants.

It has been suggested by many authors that the behavior of turbulent
two-dimensional flows should be understood by the methods of equilibrium or
non-equilibrium statistical mechanics (\cite{Onsager} \cite{Robert1} \cite%
{Someria} \cite{Caglioti} \cite{Robert2} \cite{Bouchet1} and references
therein). Modern studies in this direction concentrate in construction of
microcanonical or more general invariant Young measures, on their relation
to the small viscosity limit of the invariant measures of Navier-Stokes,
relaxation of the dynamics and phase transitions.

Here, following the inspiration provided by the results on \textit{physical
measures}, as described above, we study the stochastic stability of the
invariant measures. The plan of the paper is as follows: in Section 2,
infinitesimally invariant measures of partial differential equations are
related to the generator of the flow and in Section 3 the 2D Euler equation
with periodic boundary conditions is written as a differential equation for
its Fourier modes and it is shown that it has infinitely many invariant
measures.

In Section 4.1 we revisit the question already addressed by other authors 
\cite{Albever2} \cite{Cipriano} of whether an invariant measure of the 2D
Euler equation remains invariant when the deterministic flow is replaced by
an Ornstein-Uhlenbeck\ process. Some such measures are found, which however
correspond both to a noise perturbation and to a change of the deterministic
vector field. Therefore they are not candidates for the \textit{%
stochastically stable measures} in the sense described before. Then in
Section 4.2, we add a noise perturbation to the deterministic dynamics and
show that once a boundary condition on the dynamical space is fixed, there
is a unique measure which converges in the sense of viscosity solutions to a
measure density of the deterministic equation. This result is obtained for
the 2D Euler equation truncated to arbitrarily large $N$ Fourier modes. How
to generalize it to the infinite-dimensional case is indicated.

The result obtained in 4.2 provides a reasonable interpretation of the
stability of the large scale structures in two-dimensional fluid motion.
Because the stochastically stable invariant measure depends on the boundary
conditions (for example a cut-off at large modes), we also understand why,
depending on the particular physical environment, the structures display not
a unique but several different shapes. It also provides a plausible
explanation for the relaxation of the flow to selected structures, not as an
effect of some residual viscosity but as a result of the noise always
present in a physical system. In addition the dependence of the
stochastically stable measure on the dynamical boundary conditions might
also provide an explanation of why the same basic equation may lead to
different large scale patterns depending on the physical environment.

Finally, in Section 5, we briefly rephrase our results in configuration
space and using a recently developed stable algorithm perform a few
illustrative numerical simulations of a finite mode 2D Euler equation
perturbed by noise that show the emergence of the stochastically stable
patterns.

Most of the results in the paper refer to a truncated system, therefore to
an arbitrarily large, but finite, dynamical system. The actual extension to
an infinite system is sketched but not worked out in detail.

\section{Infinitesimally invariant measures of partial differential equations%
}

Let $\Gamma _{t}$ be the flow of a partial differential equation and $\Gamma
_{t}^{\ast }$ the push-forward semigroup acting on measures. A measure $\mu $
is invariant if%
\begin{equation}
\Gamma _{t}^{\ast }\left( \mu \right) =\mu  \label{IM1}
\end{equation}%
and infinitesimally invariant if%
\begin{equation}
\int B\varphi d\mu =0  \label{IM2}
\end{equation}%
for any differentiable function $\varphi $, $B$ being the generator of the
flow $\Gamma _{t}$. Equivalently $B^{\ast }1=0$.

Let the generator $B$ be a first or second order differential operator on a
discrete set of coordinates $\phi =\left\{ \phi _{i}\right\} $,%
\begin{equation}
B=\sum_{i,j}u_{ij}\left( \phi \right) \frac{\partial ^{2}}{\partial \phi
_{i}\partial \phi _{j}}+\sum_{i}b_{i}\left( \phi \right) \frac{\partial }{%
\partial \phi _{i}}  \label{IM3}
\end{equation}%
and consider a measure of the form\footnote{%
Here and throughout most of the paper $\prod_{i}d\omega _{i}$ stands for $%
\prod_{i=1}^{N}d\omega _{i}$ with $N$ an arbitrarily large integer. The
infinite dimensional case will be discussed in the last part of Section 4.}%
\begin{equation}
d\mu =R\left( \phi \right) \prod_{i}d\phi _{i}  \label{IM4}
\end{equation}%
To obtain the condition (\ref{IM2})%
\begin{equation*}
\int \left( B\varphi \right) R\left( \phi \right) \prod_{i}d\phi _{i}=0
\end{equation*}%
one computes the adjoint of $B$ obtaining%
\begin{eqnarray}
B^{\ast } &=&-\frac{1}{R}\left\{ \sum_{i}\frac{\partial }{\partial \phi _{i}}%
\left( Rb_{i}\right) -\sum_{i,j}\frac{\partial ^{2}}{\partial \phi
_{i}\partial \phi _{j}}\left( Ru_{ij}\right) \right\}  \notag \\
&&+\sum_{i}\left\{ -b_{i}+\frac{1}{R}\sum_{j}\frac{\partial }{\partial \phi
_{j}}\left[ R\left( u_{ij}+u_{ji}\right) \right] \right\} \frac{\partial }{%
\partial \phi _{i}}+\sum_{i,j}u_{ij}\frac{\partial ^{2}}{\partial \phi
_{i}\partial \phi _{j}}  \label{IM5}
\end{eqnarray}%
Therefore, to have $B^{\ast }1=0$, the first term in (\ref{IM5}) should
vanish leading to

\textbf{Proposition 1:} \textit{A generator }$B$\textit{\ of the form in Eq.(%
\ref{IM3}), }$u_{ij}$\textit{\ and }$b_{i}$\textit{\ being differentiable
functions, has}%
\begin{equation*}
d\mu =R\left( \phi \right) \prod_{i}d\phi _{i}
\end{equation*}%
\textit{(}$R\left( \phi \right) $\textit{\ differentiable) as an
infinitesimally invariant measure if and only if}%
\begin{equation}
\sum_{i}\frac{\partial }{\partial \phi _{i}}\left( Rb_{i}\right) -\sum_{i,j}%
\frac{\partial ^{2}}{\partial \phi _{i}\partial \phi _{j}}\left(
Ru_{ij}\right) =0  \label{IM6}
\end{equation}%
\textit{Equivalently}%
\begin{equation}
b_{i}=\frac{1}{R}\sum_{j}\frac{\partial }{\partial \omega _{j}}\left(
Ru_{ij}\right) +\frac{X_{i}}{R}  \label{IM7}
\end{equation}%
\textit{where }$X_{i}$\textit{\ is an arbitrary function satisfying }$%
\sum_{i}\frac{\partial X_{i}}{\partial \phi _{i}}=0$.

A similar result has been obtained in \cite{Airault}.

\section{The 2D Euler equation on the torus}

Consider the 2D Euler equations for an inviscid incompressible fluid 
\begin{equation}
\left\{ 
\begin{array}{l}
{\frac{\partial v}{\partial t}}=-(v\cdot \nabla )v-\nabla p \\ 
div\ v=0%
\end{array}%
\right.  \label{E1.1}
\end{equation}%
subjected to periodic boundary conditions and initial data 
\begin{equation}
v(x,0)=v_{0}(x),  \label{E1.2}
\end{equation}%
where $v(x,t)=(v_{1}(x_{1},x_{2},t),v_{2}(x_{1},x_{2},t))$ is the velocity
field of the fluid and $p=p(x,t)$ is the pressure.

Since $div\ v=0$ and $div\ v_{0}=0,$ there is a function $\psi (x,t)$ (the
stream function) such that%
\begin{equation}
v=\nabla ^{\perp }\psi =(-\partial _{x_{2}}\psi ,\partial _{x_{1}}\psi )
\label{E1.2a}
\end{equation}%
and the Euler equation becomes%
\begin{equation}
\partial _{t}\Delta \psi =-\nabla ^{\perp }\psi \cdot \nabla \Delta \psi
\label{E1.3}
\end{equation}%
As in \cite{Albever2} we consider solutions of (\ref{E1.3}) on the
2-dimensional flat torus, a square in $\mathbb{R}^{2}$ with periodic
boundary conditions, $T^{2}=[0,1]\times \lbrack 0,1]$, 
\begin{equation}
\psi (0,x_{2},t)=\psi (1,x_{2},t),\;\;\;\;\;\;\;\;\;\psi (x_{1},0,t)=\psi
(x_{1},1,t)  \label{E1.4}
\end{equation}%
$\forall x=(x_{1},x_{2})\in T^{2},\forall t\in \lbrack 0,T].$ Let us denote
by $e_{k}(x)=e^{i\,2\pi k\cdot x},$ \ $k\in \mathbb{Z}^{2}$ the
eigenfunctions for the operator $-\Delta $ with eigenvalues $4\pi ^{2}\left(
k_{1}^{2}+k_{2}^{2}\right) $, where $k\cdot x=k_{1}x_{1}+k_{2}x_{2}$. They
form a complete set of orthonormal functions in $L^{2}(T^{2})$. We expand
the solution $\psi (x,t)$ of (\ref{E1.3}) as a Fourier series%
\begin{equation*}
\psi (x,t)=\sum_{k}\phi _{k}(t)e_{k}(x).
\end{equation*}%
Since $\psi $ is a real function and we can assume $\int_{T^{2}}\psi dx=0$,
then $\phi _{-k}=\overline{\phi }_{k}$ ($\overline{z}$ being the complex
conjugate of $z$) and 
\begin{equation}
\psi (x,t)=\sum_{k\in \mathbb{Z}_{+}^{2}}\phi _{k}(t)e_{k}(x),  \label{E1.5}
\end{equation}%
where $\mathbb{Z}_{+}^{2}$ denotes the set $\{k\in \mathbb{Z}^{2}:\
k_{1}>0,\ k_{2}\in \mathbb{Z}\text{ or }k_{1}=0,\ k_{2}>0\}$.

By (\ref{E1.5}), the function $\psi $ is identified with an infinite vector
of Fourier coefficients 
\begin{equation*}
\psi =\left\{ \phi _{k}\right\} _{k\in \mathbb{Z}_{+}^{2}}
\end{equation*}%
where$\ k\in \mathbb{Z}_{+}^{2}$. We define $\mathbb{C}^{\infty }=\left\{
\phi =\left\{ \phi _{k}\right\} _{k\in \mathbb{Z}_{+}^{2}}\ :\ \phi _{k}\in 
\mathbb{C}\right\} .$

Substituting (\ref{E1.5}) in equation (\ref{E1.3}) and introducing the
operator \cite{Albever2} \cite{Albever1} \cite{Cipriano}%
\begin{equation*}
B(\phi )=\left\{ B_{k}(\phi )\right\} _{k\in \mathbb{Z}_{+}^{2}}=%
\sum_{k}B_{k}(\phi )\frac{\partial }{\partial \phi _{k}}
\end{equation*}%
with coefficients $B_{k}=B_{k}(\omega )$ 
\begin{equation}
B_{k}(\phi )=\frac{4\pi ^{2}}{k^{2}}\sum_{\substack{ h\neq k  \\ h,k\in 
\mathbb{Z}_{+}^{2}}}\left( k^{\bot }\cdot h\right) \left( k-h\right)
^{2}\phi _{h}\phi _{k-h}  \label{E1.7}
\end{equation}%
where $k^{\bot }=(-k_{2},k_{1})$, the system (\ref{E1.1}) becomes the
following infinite dimensional ordinary differential equation%
\begin{equation}
\frac{d}{dt}\phi _{k}=B_{k}(\phi )\hspace{0.5cm}k\in \mathbb{Z}_{+}^{2}
\label{E1.8}
\end{equation}%
and%
\begin{equation}
\frac{\partial B_{k}}{\partial \phi _{k}}=0  \label{E1.9}
\end{equation}

We may now find the (infinitesimally) invariant measures of the Euler
equation on the torus. For the measure (\ref{IM4}) we see from (\ref{IM3})
that with $u_{ij}\left( \phi \right) =0$, the condition (\ref{IM6}) is simply%
\begin{equation*}
\sum_{i}\frac{\partial }{\partial \phi _{i}}\left( Rb_{i}\right) =0
\end{equation*}%
that is,%
\begin{equation*}
\sum_{i}\frac{\partial }{\partial \phi _{i}}\left( RB_{i}\right) =0
\end{equation*}%
or from (\ref{E1.9})%
\begin{equation*}
\sum_{i}B_{i}\frac{\partial }{\partial \phi _{i}}R=\sum_{i}\frac{d}{dt}\phi
_{i}\frac{\partial }{\partial \phi _{i}}R=\frac{d}{dt}R=0
\end{equation*}%
In conclusion: \textit{any constant of motion of the Euler equation
generates an (infinitesimally) invariant measure.} Among them we mention the
energy $E$ and the enstrophy $S$ (or functions thereof) which in this
setting read%
\begin{equation*}
E=\frac{1}{2}\sum_{k}k^{2}\phi _{k}^{2}
\end{equation*}%
\begin{equation*}
S=\frac{1}{2}\sum_{k}k^{4}\phi _{k}^{2}
\end{equation*}%
The Poisson structure of the Euler 2D equation being degenerate, there is a
set of Casimir invariants\footnote{%
Related by Noether theorem to relabelling invariance of the fluid elements 
\cite{Salmon}} \cite{Morrison}, which are invariant for any Hamiltonian flow
with that Poisson structure. In this case they are%
\begin{equation*}
C_{f}=\int f\left( \triangle \psi \right) d^{2}x
\end{equation*}%
$f$ being an arbitrary differentiable function. Therefore there are
infinitely many invariant measures for the 2D Euler equation. The enstrophy
is the Casimir invariant for $f\left( x\right) =x^{2}$.

\section{Stochastic perturbations of the 2D Euler equation and invariant
measures}

Here we discuss stochastic stability of invariant measures in two different
settings. First, given an invariant measure of the deterministic equation,
we find the stochastic perturbation which preserves that measure when also
the deterministic part is allowed to change. Second, we discuss the
invariant measures of the stochastically perturbed system, with the
deterministic part kept fixed and also the convergence of the perturbed
measure when the perturbation tends to zero. It is this second study that is
in the spirit of the identification of the physical measure by stochastic
perturbations as it is done for finite-dimensional dissipative systems.

\subsection{Stochastic perturbations preserving a deterministic invariant
measure}

A similar such study has been performed before and we use the same setting
and notation as in \cite{Albever2} \cite{Cipriano}. We introduce the Sobolev
spaces of order $\beta \in \mathbb{R}$ on the torus $T^{2}$ 
\begin{eqnarray}
H^{\beta } &=&\left\{ \phi =\sum_{k}\phi _{k}e_{k}:\sum_{k}\left\vert
k\right\vert ^{2\beta }\left\vert \phi _{k}\right\vert ^{2}<+\infty ,\phi
_{-k}=\overset{-}{\phi _{k}}\right\}  \notag \\
&\equiv &\left\{ \phi =\left( \phi _{k}\right) _{k\in \mathbb{Z}_{+}^{2}}\in 
\mathbb{C}^{\infty }:\sum_{k\in \mathbb{Z}_{+}^{2}}\left\vert k\right\vert
^{2\beta }\left\vert \phi _{k}\right\vert ^{2}<+\infty \right\}  \label{SP0}
\end{eqnarray}%
The spaces $H^{\beta }$ are complex Hilbert spaces with inner product and
norm given by 
\begin{equation*}
<\phi ^{(1)},\phi ^{(2)}>_{H^{\beta }}=\sum_{k\in \mathbb{Z}%
_{+}^{2}}\left\vert k\right\vert ^{2\beta }\phi _{k}^{(1)}\overline{\phi }%
_{k}^{(2)},\qquad \Vert \phi \Vert _{H^{\beta }}^{2}=<\phi ,\phi >_{H^{\beta
}}.
\end{equation*}%
\textbf{Definition:} \textit{An arbitrary complex function }$f=f(\phi
):C^{\infty }\rightarrow C$\textit{\ is a cylindrical function if, for some
integer }$N$\textit{, we have }$f=f(\phi )\equiv F(\phi _{\alpha _{1}},\dots
,\phi _{\alpha _{d(N)}})$\textit{, where }$F$\textit{\ is a }$%
C_{0}^{1}(C^{N})$\textit{\ - smooth function depending only on the
components }$\phi _{\alpha _{i}}$\textit{, }$\alpha _{i}\in Z_{+,d(N)}^{2}$%
\textit{.}

Let us consider the following infinite dimensional parametric
Ornstein-Uhlenbeck operator $\varepsilon Q$ defined by%
\begin{equation}
\varepsilon Qf(\phi )=\varepsilon \sum_{k}\left\{ a_{k}\left( \phi \right) 
\frac{\partial }{\partial \phi _{k}}f(\phi )+\sigma _{k}\left( \phi \right) 
\frac{\partial ^{2}}{\partial \phi _{k}^{2}}f(\phi )\right\}  \label{SP1}
\end{equation}%
for every cylindrical function.

If we consider the operator 
\begin{equation}
Lf(\phi )=\varepsilon Qf(\phi )+\sum_{k}B_{k}(\phi )\frac{\partial }{%
\partial \phi _{k}}f(\phi )  \label{SP2}
\end{equation}%
we can see this operator as the infinitesimal generator for a stochastically
perturbed Euler flow.

Let $W(t)=\sum_{k}\frac{1}{|k|}b_{k}(t)e_{k}$ be a normalized cylindrical
brownian motion on $H^{1-\delta }$, $b_{k}(t)$ being independent copies of a
complex brownian motion. To the generator (\ref{SP2}) corresponds the
following perturbed Euler equation%
\begin{equation}
X_{k}(t)=X_{k}\left( 0\right) +\int_{0}^{t}\left\{ B_{k}\left( X\left(
s\right) \right) +\varepsilon a\left( X_{k}(s)\right) \right\}
ds+\int_{0}^{t}\sqrt{2\varepsilon \sigma _{k}\left( X_{k}(s)\right) }%
db_{k}(s),\,\,\,\,\forall k\in \mathbb{Z}_{+}^{2}.  \label{SP3}
\end{equation}

\textbf{Proposition 2:} \textit{If }$d\mu =R\left( \phi \right)
\prod_{i}d\phi _{i}$\textit{\ is an invariant measure for the (truncated)
unperturbed Euler equation, then this is also an invariant measure for the
perturbed equation (\ref{SP3}) if }$a_{k}\left( \phi \right) $\textit{\ and }%
$\sigma _{k}\left( \phi \right) $\textit{\ in (\ref{SP1}) satisfy}%
\begin{equation}
\sum_{k}\left\{ \left( a_{k}-2\frac{\partial \sigma _{k}}{\partial \phi _{k}}%
\right) \frac{\partial R}{\partial \phi _{k}}+R\left( \frac{\partial a_{k}}{%
\partial \phi _{k}}-\frac{\partial ^{2}\sigma _{k}}{\partial \phi _{k}^{2}}%
\right) -\sigma _{k}\frac{\partial ^{2}R}{\partial \phi _{k}^{2}}\right\} =0
\label{SP4}
\end{equation}

This is a direct consequence of Eq.(\ref{IM6}). As an example, for the
Gaussian measure constructed from the enstrophy%
\begin{equation}
d\mu _{S}=e^{-\frac{1}{2}\sum_{k}k^{4}\phi _{k}^{2}}\prod_{j}d\phi _{j}
\label{SP5}
\end{equation}%
Eq.(\ref{SP4}) is satisfied by%
\begin{equation}
a_{k}=-k^{2}\phi _{k},\hspace{0.5cm}\sigma _{k}=\frac{1}{k^{2}}  \label{SP6}
\end{equation}%
and for the Gaussian measure constructed from the renormalized energy%
\begin{equation}
d\mu _{E}=e^{-:E:}\prod_{j}d\phi _{j}  \label{SP7}
\end{equation}%
\begin{equation}
a_{k}=-\phi _{k},\hspace{0.5cm}\sigma _{k}=\frac{1}{k^{2}}  \label{SP8}
\end{equation}%
where $:E:=\frac{1}{2}\left( \sum_{k}k^{2}\phi _{k}^{2}-\mathbb{E}\left[
\sum_{k}k^{2}\phi _{k}^{2}\right] \right) $.

Notice that in (\ref{SP5}) and (\ref{SP7}) we are considering a truncation
of the 2D Euler equation to arbitrarily large $N$ modes. In the $%
N\rightarrow \infty $ limit the flat measure $\prod_{j}d\phi _{j}$ makes no
sense and another reference measure should be used.

One sees that for these invariant measures of the unperturbed Euler
equation, there are specific Ornstein-Uhlenbeck perturbations that preserve
it as an invariant measure. However, in each case we are not only adding
noise but also modifying the deterministic part. In the first (enstrophy)
case we are actually adding noise to a Navier-Stokes equation%
\begin{equation*}
\begin{array}{ccc}
\partial _{t}\Delta \psi & = & -\nabla ^{\perp }\psi \cdot \nabla \Delta
\psi +\varepsilon \triangle ^{2}\psi \\ 
{\frac{\partial v}{\partial t}} & = & -(v\cdot \nabla )v+\varepsilon
\triangle v-\nabla p%
\end{array}%
\end{equation*}%
and in the renormalized energy case%
\begin{equation*}
\begin{array}{ccc}
\partial _{t}\Delta \psi & = & -\nabla ^{\perp }\psi \cdot \nabla \Delta
\psi -\varepsilon \triangle \psi \\ 
{\frac{\partial v}{\partial t}} & = & -(v\cdot \nabla )v-\varepsilon
v-\nabla p%
\end{array}%
\end{equation*}%
Therefore, because invariance of these measures requires a fine tuning with
both the deterministic and the stochastic components being modified with the
same intensity $\varepsilon $, they do not seem to be the right candidates
for the physical measures of the 2D Euler equation. The same applies to the
results of Kuksin \cite{Kuksin} who, using a viscosity of intensity $%
\varepsilon $ and a $\sqrt{\varepsilon }$ noise, shows that the collection
of unique invariant measures so obtained is tight and converges in the $%
\varepsilon \rightarrow 0$ limit to a measure of the deterministic Euler
equation.

Incidentally, also the microcanonical measures, that have been studied by a
number of authors, do not seem to qualify as stochastically stable measures
even with reasonable modifications of the deterministic part of the equation.

That the selection of a unique invariant measure requires a fine tuning, of
both the noise and the deterministic terms, makes these, otherwise
interesting, results irrelevant for the interpretation of physical
phenomena, where such fine tuning is not to be expected.

\subsection{The zero noise limit of the invariant measure of a stochastic
system}

In the previous subsection we have dealt with stochastic perturbations which
preserve invariant measures of (\ref{E1.8}). As stated before, of more
interest for the characterization of the \textit{physical measures} would be
to find noise-perturbed systems with an unique invariant measure and to
construct the zero-noise limit of that measure. This we discuss now, not for
the infinite dimensional system but again for its Galerkin approximations of
arbitrary order $N$ \cite{Albever3}%
\begin{equation}
\frac{d}{dt}\phi _{k}=B_{k}^{N}(\phi )\hspace{0.5cm}k\in \mathbb{Z}_{0}^{2}%
\hspace{0.8cm}\left\vert k\right\vert \leq N  \label{V1}
\end{equation}%
\begin{equation}
B_{k}^{N}(\phi )=\frac{4\pi ^{2}}{k^{2}}\sum_{\substack{ 0<\left\vert
h\right\vert \leq N  \\ 0<\left\vert k-h\right\vert \leq N}}\left( k^{\bot
}\cdot h\right) \left( k-h\right) ^{2}\phi _{h}\phi _{k-h}  \label{V2}
\end{equation}%
When noise is added to (\ref{V1}), without changing the deterministic part,
the equation for the density $R\left( \phi \right) $ of the invariant
measure becomes%
\begin{equation}
\sum_{k}B_{k}^{N}\left( \phi \right) \frac{\partial }{\partial \phi _{k}}%
R-\varepsilon \sigma _{k}\frac{\partial ^{2}}{\partial \phi _{k}^{2}}R=0
\label{VHJ}
\end{equation}%
Two cases are of physical interest, namely $\sigma _{k}=1$ and $\sigma _{k}=%
\frac{1}{k^{2}}$, corresponding respectively to a uniform noise in all
Fourier modes or to a decreasing noise intensity in higher modes. However,
by the change of variables $z_{k}=\left\vert k\right\vert \phi _{k}$ and $%
B_{k}^{N^{\prime }}\left( \phi \right) =\left\vert k\right\vert
B_{k}^{N}\left( \phi \right) $ the second case becomes identical to the
first one and we have to deal with%
\begin{equation}
\sum_{k}B_{k}^{N^{\prime }}\left( z\right) \frac{\partial }{\partial z_{k}}%
R-\varepsilon \frac{\partial ^{2}}{\partial z_{k}^{2}}R=0  \label{V3}
\end{equation}%
which we recognize as an elliptic regularization of a first order
Hamilton-Jacobi equation. As shown before, this Hamilton-Jacobi equation $%
\left( \varepsilon =0\right) $ has at least as many generalized solutions as
the number of constants of motion of the $N-$Galerkin approximation to the
Euler equation. Hence, existence and uniqueness of a stochastically-stable
solution for $R$ is equivalent to the establishment of a viscosity solution%
\footnote{%
A viscosity solution is a weak solution which need not be everywhere
differentiable (see \cite{Lions1}).} for this Hamilton-Jacobi problem \cite%
{Lions1} \cite{Lions2} \cite{Lions3}, in particular in its vanishing
viscosity modality \cite{Lions3} \cite{Evans} (ch. 10).

However, the solution of this problem is strongly depend on the domain where
the $R$ function is defined, therefore on the dynamical boundary conditions.
What this means in practical terms is that the fluid under study might not
be exploring all possible intensities in all modes. In Eq.(\ref{V3}) this
would be coded by particular boundary conditions on the $R$ function.

Associated to the uniformly elliptic equation (\ref{V3}) there is a
diffusion process $X_{\varepsilon }\left( t\right) $ with diffusion
coefficient $\sqrt{\varepsilon }$ and drift $B_{k}^{N^{\prime }}\left(
z\right) $. In each bounded domain $D$ of $z-$space, the drift, being a
quadratic polynomial, is uniformly Lipschitz continuous. Therefore the
Dirichlet problem of Eq.(\ref{V3}) has a unique solution with stochastic
representation%
\begin{equation}
\left. R_{\varepsilon }\left( z\right) \right\vert _{D}=\mathbb{E}%
_{z}\left\{ f\left( X_{\varepsilon }\left( \tau \right) \right) \right\}
\label{V4}
\end{equation}%
$f$ being the boundary condition at $\partial D$ and $\tau $ the first exit
time from $D$ (\cite{Friedman} ch. 6).

For a bounded smooth boundary condition the solution $R_{\varepsilon }$ in (%
\ref{V4}) is bounded and continuous on compact subsets of $D$. Then, when $%
\varepsilon \rightarrow 0$ $R_{\varepsilon }$ converges locally uniformly to
a function $R$. This function is not necessarily a classical solution of $%
\sum_{k}B_{k}^{N^{\prime }}\left( z\right) \frac{\partial }{\partial z_{k}}%
R=0$, but a standard construction (\cite{Evans}, ch.10) shows that it is a
viscosity solution, in the sense that, given a $\mathbb{C}^{\infty }$
function $g$, if $R-g$ has a local maximum at a point $z_{0}$ then $%
\sum_{k}B_{k}^{N}\left( z_{0}\right) \frac{\partial }{\partial z_{k}}g\left(
z_{0}\right) \leq 0$ and if it is a local minimum $\sum_{k}B_{k}^{N}\left(
z_{0}\right) \frac{\partial }{\partial z_{k}}g\left( z_{0}\right) \geq 0$.
Hence,

\textbf{Proposition 3:} \textit{For each choice of boundary conditions in }$%
z-$\textit{\ space and noise level (}$\varepsilon $\textit{), one has a
unique measure density }$R_{\varepsilon }\left( z\right) $\textit{, solution
of (\ref{V3}). Furthermore, in the }$\varepsilon \rightarrow 0$\textit{\
limit, }$R_{\varepsilon }$\textit{\ converges to a viscosity solution of }$%
\sum_{k}B_{k}^{N^{\prime }}\left( z\right) \frac{\partial }{\partial z_{k}}%
R=0$.

For consistency with the $\varepsilon =0$ case, it is convenient to have the
boundary function at each $\partial D_{n}$ constructed from a constant of
motion of the 2D Euler equation, for example the enstrophy ($\left.
f_{n}\right\vert _{\partial D_{n}}=e^{-\frac{1}{2}\sum_{k}k^{4}\phi
_{k}^{2}} $) as in (\ref{SP5}). Then the viscosity solution would provide a
measure density which for very large mode amplitudes behaves like the
enstrophy measure. In this construction the measures may be made to coincide
in the boundary with one of the infinitely many invariant measures discussed
in section 2. However in the interior of the specified domain the
stochastically stable solution will not in general coincide with the
solution chosen for the boundary. Also, the solution that is obtained is not
in a strict sense an invariant measure for the original equation because of
the limitations put on the domain by the boundary conditions. However it
follows from (\ref{V4}) that, for a positive boundary condition, $R$ is a
positive density.

So far we have dealt with $N$-dimensional Galerkin approximations to the 2D
Euler equation. When $N\rightarrow \infty $ several modifications are
needed. The first one is in the equation (\ref{IM4}) because it makes no
sense to define $R\left( \phi \right) $ as a density of the non-existent
flat measure in infinite dimensions. Instead, $R\left( \phi \right) $ should
be defined as the Radon-Nykodim derivative for some other measure, for
example the Gaussian enstrophy measure. Then the equation for the density $%
R\left( \phi \right) $ would be%
\begin{equation}
\sum_{k}\left\{ B_{k}\left( \phi \right) \frac{\partial }{\partial \phi _{k}}%
-k^{4}\phi _{k}B_{k}\left( \phi \right) \right\} R\left( \phi \right) =0
\label{V5}
\end{equation}%
an Hamilton-Jacobi equation in infinite dimensions. Such equations have been
extensively studied \cite{Lions-infinite} and given the appropriate boundary
condition, for example $R\left( \phi \right) \rightarrow 1$ for large $%
\left\vert \phi \right\vert $, the construction of the density as a limiting
viscosity solution of%
\begin{equation}
\sum_{k}\left\{ B_{k}\left( \phi \right) \frac{\partial }{\partial \phi _{k}}%
-k^{4}\phi _{k}B_{k}\left( \phi \right) -\varepsilon \frac{\partial ^{2}}{%
\partial \phi _{k}^{2}}\right\} R\left( \phi \right) =0  \label{V6}
\end{equation}%
would follow similar steps as in the finite dimensional case.

Proposition 3 establishes the existence of stochastically stable measures as
viscous solutions of an elliptic regularized Hamilton-Jacobi equation.. The
solutions are defined once the boundary conditions at large $\phi
_{k}^{\prime }s$ are fixed, for example, by some invariant measure of the
deterministic 2D Euler equation.

In conclusion, the present result provides an interpretation of the
stability of the large coherent structures in two dimensional fluid motion
somewhat different from what has been suggested in the past. Some past
treatments start from the fact that the stationary points of constants of
motion are steady state solutions and choose an appropriate linear
combination $G$ of the constants of motion as a potential and adding to the
equations a $-\alpha G$ term develop a dissipative Langevin dynamics.
Alternatively, other approaches look for maxima of the entropy, which of
course depend on a previous choice of measure. In particular the
microcanonical measure, that has been favored, is not a solution of the
elliptic regularization of the Hamilton-Jacobi equation for finite noise
level $\varepsilon $. Whether it can, in some sense, be identified with a
viscosity solution in the $\varepsilon \rightarrow 0$ limit is an open
question.

In contrast with previous interpretations, our analysis suggests that the
coherent structures observed in 2D hydrodynamics are associated to
configurations that are stochastically stable measures uniquely determined
by the boundary conditions in $\left\{ \phi \right\} -$ space. Some authors
have suggested that the convergence of two-dimensional fluid dynamics to
stable or quasi-stable large scale structures is associated to dissipative
effects. Of course, a dissipative effect may be interpreted as a dynamical
boundary condition, for example a suppression of the high Fourier modes. But
what our result shows is that uniqueness of the invariant measure is
associated to the dynamical boundary conditions, dissipative or otherwise.

\section{Stochastically stable configurations: Numerical illustrations}

Here, instead of the Fourier mode decomposition and truncation we use
configuration space variables. Corresponding to the Fourier mode truncation,
one has the stream function defined at a grid of $N\times N$ points.
Therefore instead of Fourier modes, one has values of the stream function at
points in a grid and the same type of results are expected. The truncated
equation is%
\begin{equation}
\partial _{t}\left( \Delta \psi \right) _{ij}=-\left( \nabla ^{\perp }\psi
\cdot \nabla \right) _{ik}\left( \Delta \psi \right) _{kj}  \label{5.1}
\end{equation}%
where now $\Delta $ and $\nabla $ stand for the discrete Laplacian and
discrete gradient. The evolution of the stream function is obtained by the
inversion of a Poisson equation%
\begin{equation}
\psi _{ij}=\left( \Delta ^{-1}\right) _{ik}\left( \Delta \psi \right) _{kj}
\label{5.2}
\end{equation}%
with the physically irrelevant condition%
\begin{equation}
\sum_{ij}\psi _{ij}=0  \label{5.3}
\end{equation}%
What has been proved in the previous section was the existence of unique
stochastically stable \textit{measures} once the dynamical boundary
conditions are fixed, not the existence of unique stochastically stable 
\textit{solutions}. However it is to be expected that, when perturbed by
small noise, the solutions will be concentrated on the regions where the
measure is maximal. This is now illustrated with numerical simulations. To
perform these simulations in a reliable way one should insure that the
observed effects come from the noise perturbations and not from round-off or
numerical instabilities of the algorithm. In this case the evolution
operator $M$%
\begin{equation*}
M=\nabla ^{\perp }\psi \cdot \nabla
\end{equation*}%
a $N^{2}\times N^{2}$ matrix, is problematic because for general values of $%
\psi $ it may have both singular values greater and smaller than one.
Therefore neither an explicit nor an implicit scheme would be stable. The
solution is found by splitting $M$ into 
\begin{equation*}
M=M_{1}+M_{2}
\end{equation*}%
in such a way that the singular values of both $\left( 1-M_{1}\right) $ and $%
\left( 1+M_{2}\right) ^{-1}$ are $\leq 1$. This provides a semi-implicit
scheme \cite{Bizarro} which is stable or marginally stable.

\begin{figure}[htb]
\centering
\includegraphics[width=0.5\textwidth]{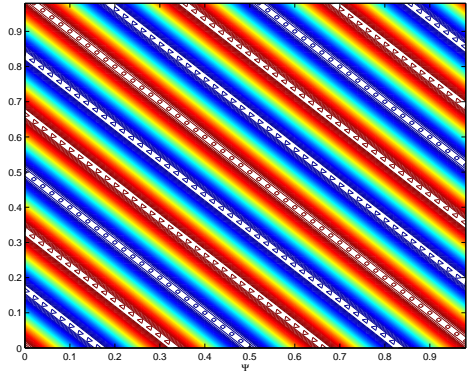}
\caption{Initial condition: a pure Fourier mode}
\label{initial_mode}
\end{figure}

\begin{figure}[htb]
\centering
\includegraphics[width=0.5\textwidth]{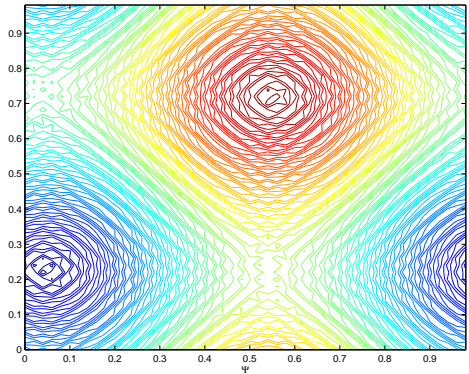}
\caption{Pattern obtained from the one in Fig.1 after evolution with noise}
\label{condensation}
\end{figure}

The semi-implicit algorithm was used with initial condition corresponding to
a single Fourier mode (Fig.\ref{initial_mode}), which is a stationary
solution of (\ref{5.1}-\ref{5.3}). However when noise is added, the solution
becomes unstable and converges to an almost stable pattern as shown in Fig.%
\ref{condensation}.

One sees that the pattern is close to the density of the first Fourier mode.
The configuration is not unique. For different runs of the simulation one
obtains essentially the same pattern but in different positions on the
torus, always close to a first Fourier mode with different phases. This
condensation in the first mode, first observed by Kraichnan and Montgomery 
\cite{Kraichnan}, has been discussed before in the framework of a
energy-enstrophy microcanonical measure \cite{Bouchet2}. However, although
we are in a finite $N$ setting, no hint of the microcanonical distribution
is apparent. For this first simulation no limitation is put on the dynamical
variable, meaning that the dynamical space is $\mathbb{R}^{N^{2}}$. Unique
solutions of the measure equation (\ref{V3}) of the type (\ref{V4}) do not
apply. However uniqueness of the solution in the $\mathbb{R}^{N^{2}}$ case
are also to be expected \cite{Albever4}.

\begin{figure}[htb]
\centering
\includegraphics[width=0.5\textwidth]{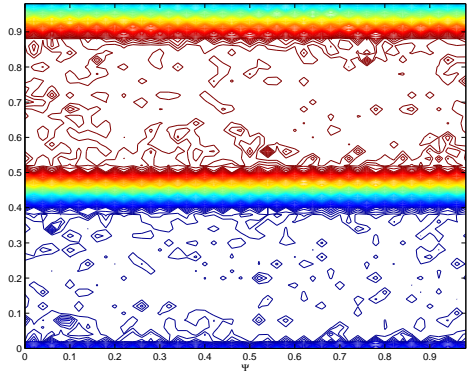}
\caption{Pattern obtained when the stream function magnitude has an upper
bound}
\label{truncated}
\end{figure}

\begin{figure}[htb]
\centering
\includegraphics[width=0.5\textwidth]{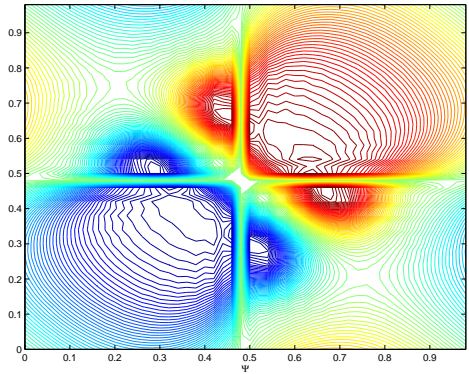}
\caption{Pattern obtained when the stream function is pinned down to zero at
two lines}
\label{cross}
\end{figure}

To explore different boundary conditions in the dynamical space, we
considered a case where the values of the stream functions are constrained
to be in a box and a case where the stream function is constrained to be
zero along two orthogonal lines. We started again from a large mode solution
which evolves under noise. The results are shown in Figs.\ref{truncated} and %
\ref{cross}. Notice that for simplicity we have considered boundary
conditions on the stream function, not on physical velocities which are
related to the stream function by Eq.(E1.2a). Boundary conditions on the
physical velocities would correspond to boundary conditions on the
derivatives of the stream function.

In this paper we have argued for the relevance of stochastically stable
measures as the generators of the coherent structures observed in (quasi)
two dimensional fluid flows. However most of our results are based on
Galerkin approximations of arbitrary but nevertheless finite dimension. In
spite of the intuition provided by Eq.(\ref{V6}), the infinite dimension
limit characterization remains, of course, an open question.

An alternative approach to the establishment of invariant measures in 2D
fluid dynamics has been the Young measure and point vertex model with finite
or variable number of vortices \cite{Robert2} \cite{Kraichnan} \cite%
{Benfatto} \cite{Sire} \cite{Xavier}, which goes back to the pioneering work
of Onsager \cite{Onsager}. In this approach, where infinite $N$ limits have
been established, Gibbs measures of the vortex model may be identified with
coherent structures, however, the selection role of stochastic stability to
choose among a basically infinite set of measures is not so clear.


\begin{thebibliography}{99}
\bibitem{Sinai} Ya. G. Sinai; \textit{Gibbs measures in ergodic theory},
Russian Math. Surveys 27 (1972) 21-69.

\bibitem{Bowen} R. Bowen; \textit{Equilibrium states and the ergodic theory
of Anosov diffeomorphisms}, Springer Lecture Notes in Math. 470 (1975)

\bibitem{Ruelle1} D. Ruelle; \textit{A measure associated with Axiom A
attractors}, Amer. J. Math. 98 (1976) 619-654.

\bibitem{Ruelle2} D. Ruelle;\textit{\ Chaotic evolution and strange
attractors}, Cambridge University Press 1989.

\bibitem{Katok} A. Katok and B. Hasselblatt; \textit{Introduction to the
modern theory of dynamical systems}, Cambridge University Press 1995.

\bibitem{Barreira} L. Barreira and Y. Pesin; \textit{Smooth ergodic theory
and nonuniformly hyperbolic dynamics}, in "Handbook of dynamical systems 1B,
chapter 2, pp. 57-264", B. Hasselblatt and A. Katok (Eds.), Elsevier,
Amsterdam 2006.

\bibitem{Kifer1} Y. I. Kifer; \textit{On small random perturbations of some
smooth dynamical systems}, Mat. URSS Izv. 8 (1974) 1083-1107.

\bibitem{Young} L.-S. Young; \textit{Stochastic stability of hyperbolic
attractors}, Ergodic Theory Dynam. Systems 6 (1986) 311--319.

\bibitem{Kifer2} Y. I. Kifer; \textit{Random Perturbations of Dynamical
Systems}, Birkh\"{a}user 1988.

\bibitem{Viana1} J. F. Alves and M. Viana; \textit{Statistical stability for
robust classes of maps with non-uniform expansion}, Ergodic Theory and
Dynamical Systems 22 (2002) 1-32.

\bibitem{Viana2} M. Benedicks and M. Viana; \textit{Random perturbations and
statistical properties of H\'{e}non-like maps}, Ann. Inst. H. Poincar\'{e}
Anal. Non Lin\'{e}aire 23 (2006) 713--752.

\bibitem{Alves} J. F. Alves, M. Carvalho and J. M. Freitas; \textit{%
Statistical stability for H\'{e}non maps of the Benedicks--Carleson type},
Ann. Inst. H. Poincar\'{e} Anal. Non Lin\'{e}aire 27 (2010) 595--637.

\bibitem{Bouchet1} F. Bouchet and A. Venaille; \textit{Statistical mechanics
of two-dimensional and geophysical flows}, Physics Reports 515 (2012)
227-295.

\bibitem{Bouchet2} F. Bouchet and M. Corvellec; \textit{Invariant measures
of the 2D Euler and Vlasov equations}, J. Stat.\ Mechanics: Theory and
Experiment (2010) P08021.

\bibitem{Onsager} L. Onsager; \textit{Statistical hydrodynamics}, Nuovo
Cimento supl. 6 (1949) 249-286.

\bibitem{Robert1} R. Robert; \textit{A maximum entropy principle for
two-dimensional perfect fluid dynamics}, J. Stat. Phys. 65 (1991) 531-553.

\bibitem{Someria} R. Robert and J. Sommeria; \textit{Statistical equilibrium
states for two-dimensional flows}, J. Fluid Mech. 229 (1991) 291-310.

\bibitem{Caglioti} E. Caglioti, P. L. Lions, C. Marchioro and M. Pulvirenti; 
\textit{A special class of stationary flows for two-dimensional Euler
equations: A statistical mechanics description}, Commun. Math. Phys. 174
(1995) 229-260.

\bibitem{Robert2} R. Robert; \textit{On the statistical mechanics of 2D
Euler equation}, Commun. Math. Phys. 212 (2000) 245-256.

\bibitem{Albever2} S. Albeverio and A. B. Cruzeiro, \textit{Global flows
with invariant (Gibbs) measures for Euler and Navier-Stokes two dimensional
fluids}, Commun. Math. Phys. 129 (1990) 431-444.

\bibitem{Cipriano} F. Cipriano, \textit{The two dimensional Euler equation:
a statistical study}, Commun. Math. Phys. (1999) 139-154.

\bibitem{Airault} H. Airault and H. Ouerdiane; \textit{Invariant measure for
some differential operators and unitarizing measure for the representation
of a Lie group. Examples in finite dimension}, Banach Center Publications 96
(2012) 11-34.

\bibitem{Albever1} S.Albeverio, M. Ribeiro de Faria and R. Hoegh-Krohn, 
\textit{Stationary measures for the periodic Euler flow in two dimensions},
J. Stat. Phys. 20 (1979) 585-595.

\bibitem{Morrison} P. J. Morrison; \textit{Hamiltonian description of the
ideal fluid}, Rev. Modern Phys. 70 (1998) 467-521.

\bibitem{Salmon} R. Salmon; \textit{Lectures on Geophysical Fluid Dynamics},
Oxford Univ. Press 1998.

\bibitem{Kuksin} S. B. Kuksin; \textit{The Eulerian limit for 2D statistical
hydrodynamics}, J. Statist. Phys. 115 (2004) 469--492.

\bibitem{Albever3} S. Albeverio and B. Ferrario; \textit{Some methods of
infinite dimensional analysis in hydrodynamics: An introduction}, in "SPDE
in Hydrodynamic: Recent Progress and Prospects", G. Da Prato and M. R\"{o}%
ckner (Eds.), Springer, Berlin 2008.

\bibitem{Lions1} M. G. Crandall and P. L.. Lions; \textit{Viscosity
solutions of Hamilton-Jacobi equations}, Trans. Amer. Math. Soc. 277 (1983)
1-42.

\bibitem{Lions2} M. G. Crandall, L. C. Evans and P. L. Lions; \textit{Some
properties of viscosity solutions of Hamilton-Jacobi equations}, Trans.
Amer. Math. Soc. 282 (1984) 487-502.

\bibitem{Lions3} P. L. Lions; \textit{Generalized solutions of
Hamilton-Jacobi equations}, Pitman, London 1982.

\bibitem{Evans} L. C. Evans; \textit{Partial differential equations},
American Mathematical Society, Providence, R. I. 2010.

\bibitem{Albever4} S. Albeverio, V. Bogachev and M. R\"{o}ckner; \textit{On
Uniqueness of Invariant Measures for Finite- and Infinite-Dimensional
Diffusions}, Comm. on Pure and Applied Math. 52 (1999) 325--362.

\bibitem{Friedman} A. Friedman; \textit{Stochastic differential equations
and applications}, vol. 1, Academic Press, New York 1975.

\bibitem{Lions-infinite} M. G. Crandall and P. L. Lions; \textit{%
Hamilton-Jacobi equations in infinite dimensions, }J. Funct. Anal. 62
(1985), 379-396; 65 (1986), 368-405; 68 (1986) 214-247; 90 (1990) 237-283;
97 (1991) 417-465; 125 (1994) 111-148.

\bibitem{Bizarro} J. P. Bizarro, L. Ven\^{a}ncio and R. Vilela Mendes; 
\textit{A stable semi-implicit algorithm}, arXiv:1905.04520.

\bibitem{Kraichnan} R. Kraichnan and D. Montgomery; \textit{Two-dimensional
turbulence}, Rep. Prog. Phys. 43(1980) 547-619.

\bibitem{Benfatto} G. Benfatto, P. Picco and M. Pulvirenti; \textit{On the
invariant measures for the two-dimensional Euler flow}, J. Statistical
Physics 46 (1987) 729-742.

\bibitem{Sire} C. Sire and P.-H. Chavanis; \textit{Numerical renormalization
group of vortex aggregation in two-dimensional decaying turbulence: The role
of three-body interactions}, Phys. Rev. E 61 (2000) 6644-6653.

\bibitem{Xavier} X. Leoncini, A. Barrat, C. Josserand and S.
Villain-Guillot; \textit{Offsprings of a point vortex}, Eur. Phys. J. B 82
(2011) 173-178.
\end{thebibliography}
\end{document}